\documentclass[aos,MSNbibl,nameyear,rotating,dvips]{arximspdf}
\usepackage{graphicx}
\doi{10.1214/12-AOS980} 
\referstodoi{10.1214/11-AOS949}
\volume{40}
\issue{4}
\pubyear{2012}
\firstpage{1973}
\lastpage{1977}

\makeatletter
\renewcommand{\citep}[1]{[\citeauthor{#1} (\citeyear{#1})]}
\newcommand{\citepp}[1]{\citeauthor{#1} (\citeyear{#1})}
\makeatother

\begin{document}
\begin{frontmatter}

\title{Discussion: Latent variable graphical model selection
via convex optimization}

\runtitle{Comment}

\begin{aug}
\author[A]{\fnms{Steffen} \snm{Lauritzen}}
\and
\author[A]{\fnms{Nicolai} \snm{Meinshausen}\corref{}\ead[label=e2]{meinshausen@stats.ox.ac.uk}}

 \runauthor{S. Lauritzen and N. Meinshausen}
 \affiliation{University of Oxford}
\address[A]{Department of Statistics\\
University of Oxford\\
1 South Parks Road\\
Oxford, OX1 3TG\\
United Kingdom\\
\printead{e2}} 
\end{aug}

 \received{\smonth{2} \syear{2012}}



\end{frontmatter}

We want to congratulate the authors for a thought-provoking and very
interesting paper. Sparse modeling of the concentration matrix has
enjoyed popularity in recent years. It has been framed as a
computationally convenient convex $\ell_1$-constrained estimation
problem in \citet{yuan05model} and can be applied readily to
higher-dimensional problems. The authors argue---we think correctly---that the sparsity of the concentration matrix is for many applications
more plausible after the effects of a few latent variables have been
removed. The most attractive point about their method is surely that it
is formulated as a convex optimization problem. Latent variable fitting
and sparse graphical modeling of the conditional distribution of the
observed variables can then be obtained through a single fitting
procedure.

\section*{Practical aspects}

The method deserves wide adoption, but this will only be realistic if
software is made available, for example, as an R-package. Not many
users will go to the trouble of implementing the method on their own,
so we will strongly urge the authors to do so.

\section*{An imputation method}

In the absence of readily available software, it is worth thinking
whether the proposed fitting procedure can be approximated by methods
involving known and well-tested computational techniques. The
concentration matrix of observed and hidden variables is
\[
 K =\pmatrix{  K_O & K_{OH} \cr K_{HO} &  K_H  },
\]
where we have deviated from the notation in the paper by omitting the
asterisk. The proposed estimator $\hat{S}_n=\hat{K}_O$  of $K_O$ was
defined as
\begin{eqnarray} \label{eq:orig}
&&(\hat{K}_O,\hat{L}_n) =  \operatorname{argmin}_{S,L}   -\ell(S-L;\Sigma_O^n)   +  \lambda_n\bigl(\gamma\|S\|_1+\operatorname{tr}(L)\bigr)  \\
&&\hspace*{191pt}\mbox{such that  } S-L \succ 0, L\succ
0,\hspace*{-191pt}
\end{eqnarray}
where $\Sigma^n_O$ is the empirical covariance matrix of the observed
variables.\vadjust{\goodbreak}

An alternative would be to replace the nuclear-norm penalization with a
fixed constraint $\kappa$ on the rank of the hidden variables,
replacing problem (\ref{eq:orig}) with
\begin{eqnarray}\label{eq:margpenlikelihood}
&&(\hat{K}_O,\hat{L}_n )  =  \operatorname{argmin}_{S,L} -\ell (S-L; \Sigma^n_O)  +  \lambda_n \|S\|_1 \nonumber  \\[-8pt]\\[-8pt]
&& \eqntext{\mbox{such that } S-L \succ 0 \mbox{ and } L\succ 0 \mbox{
and }\operatorname{rank}(L)\le \kappa .}
\end{eqnarray}

This can be achieved by a missing-value formulation in combination with
use of the EM algorithm, which also applies in a penalized likelihood
setting \citep{green1990use}. Let the hidden variables be of a fixed
dimensionality $\kappa$ and assume for a moment these are observed so
one would find the concentration matrix $\hat{K}$ of the joint
distribution of the observed variables $X_O$ and hidden variables $X_H$
based on the complete data penalized likelihood as
\begin{equation} \label{eq:penlikelihood}
\operatorname{argmin}_{K}   -\log f_K(X_O,X_H) + \lambda \| K_O\|_1 ,
\end{equation}
where $f_K$ is the joint density of $(X_O,X_H)$. This formulation is
very similar to the missing-value problem treated in
\citet{stadler2012missing}, except for the fact that we only penalize
the concentration matrix $K_O$ of the observed variables, in analogy
with the proposed latent-variable approach. The  EM algorithm
iteratively  replaces the  likelihood in (\ref{eq:penlikelihood}) for
$t=1,\ldots,T$ by its conditional expectation and thus finds
$\hat{K}^{t+1}$ as
\begin{equation}\label{eq:iter}
\hat{K}^{t+1}  =  \operatorname{argmin}_K  - E _{\hat{K}^{t}} \{ \log
f_K(X_O,X_H) | X_O \} + \lambda \|K_O\|_1 .
\end{equation}
The iteration is guaranteed not to increase the negative marginal
penalized likelihood at every stage and will, save for
unidentifiability, converge to the minimizer in
(\ref{eq:margpenlikelihood}) for most starting values. Without loss of
generality, one can fix the conditional concentration matrix $K_H$ of
the hidden variables to be the identity so that these are conditionally
independent with variance $1$, given the observed variables. Then
$-K_{OH}$ is equal to the  regression coefficients of the observed
variables on the hidden variables. As starting value we have let
$-\hat{K}^{0}_{OH}$ be equal to these with hidden variables determined
by a principal component analysis.

The expectation in (\ref{eq:iter}) can be written as the log-likelihood
of a Gaussian distribution with concentration matrix $K$ and empirical
covariance matrix~$W^t$, where
\[
 W^t = \pmatrix{  \Sigma^n_O &  -  \Sigma^n_O \hat K^t_{OH} \cr   -  \hat K^t_{HO}\Sigma^n_O &  \mathbf I +
\hat K^t_{HO} \Sigma^n_O \hat K^t_{OH}}.
\]
The
sufficient statistics involving the missing data are  thus ``imputed''
in~$W^t$. Each of the updates (\ref{eq:iter}) can now be computed with
the \emph{graphical lasso} \citep{friedman2007sic}.

\begin{figure}

\includegraphics{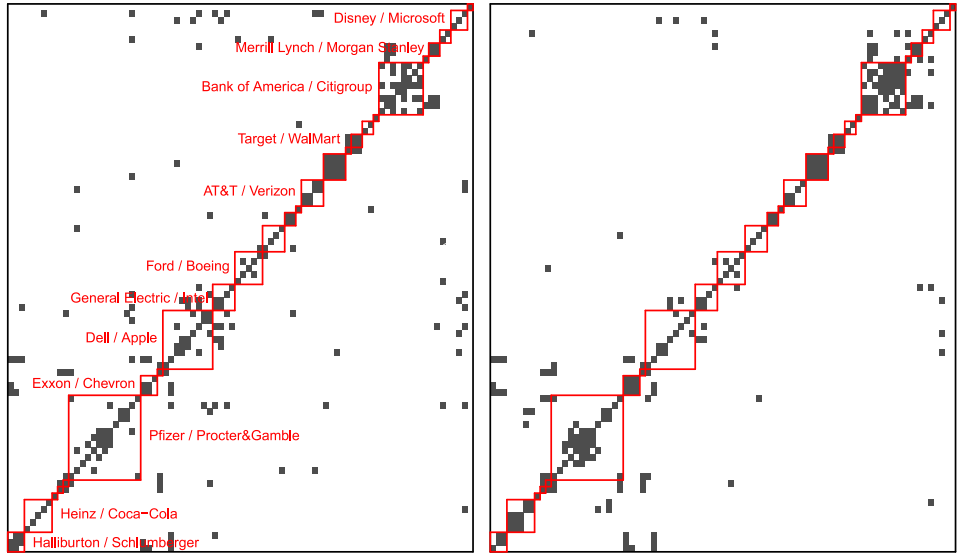}
 \caption{The nonzero entries of the concentration
matrix $\hat{K}_O$, using the proposed procedure~(\protect\ref{eq:orig})
(left) and the imputation method in
(\protect\ref{eq:penlikelihood}) (right). Two representative companies
are shown for some of the sectors.}\label{fig:conc}\vspace*{-3pt}
\end{figure}

We thought it would be interesting to compare the two methods on the
data example given in the paper. Figure \ref{fig:conc} shows the
solution $\hat{K}_O$ for the stock-return example when using the
proposed method (\ref{eq:orig}) and the imputation method
(\ref{eq:penlikelihood}) with 4 iterations. The number $\kappa$ of
latent variables\vadjust{\goodbreak} and the number of nonzero edges in $\hat{K}_O$ is
adjusted to be the same as in the original estimator.

\begin{sidewaysfigure}

\includegraphics{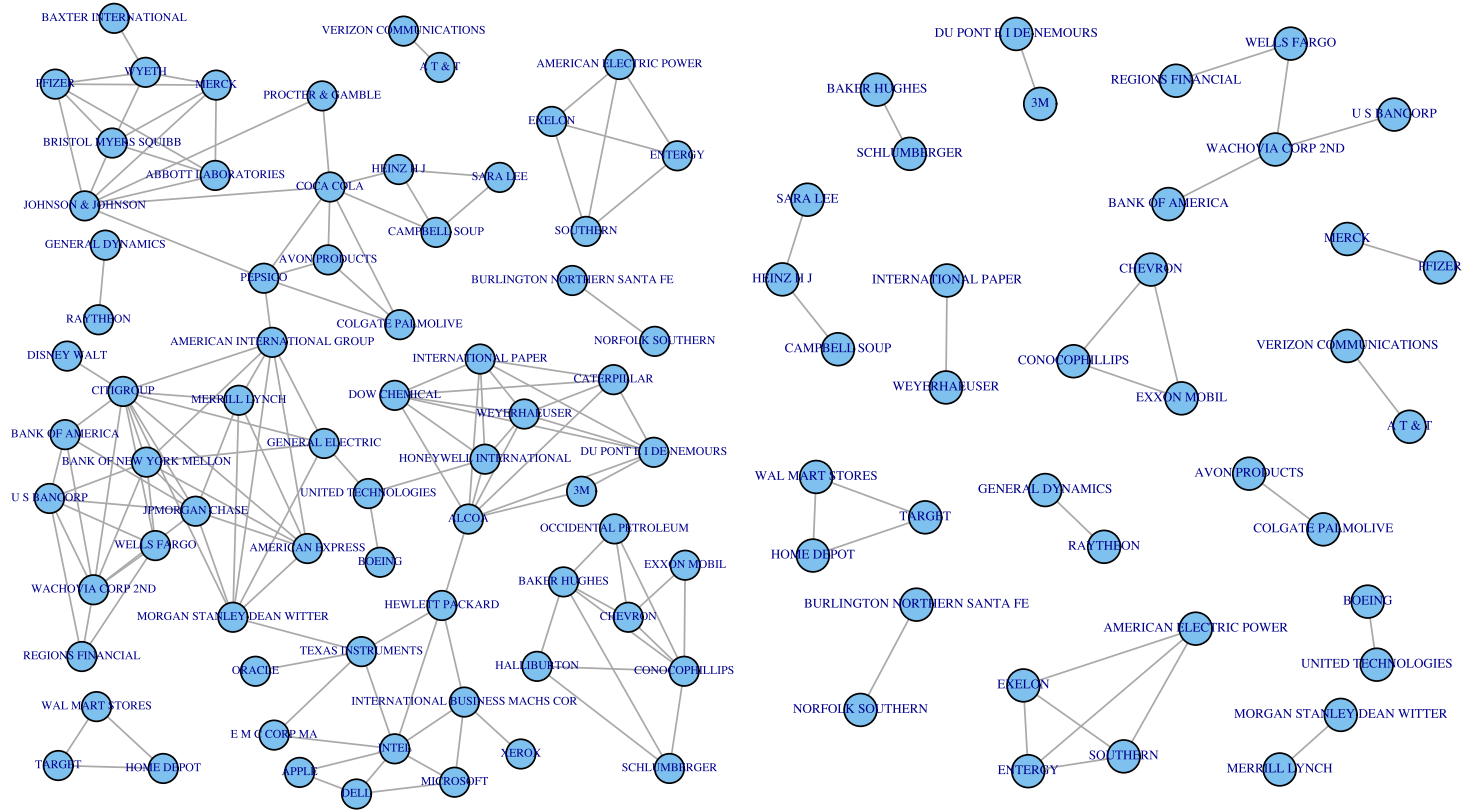}
 \caption{Left: the graph of the imputation method as
in (\protect\ref{eq:penlikelihood}). Right: the graph of the stable
edges. In both cases, isolated vertices have been removed from the
display.}\label{fig:stable}
\end{sidewaysfigure}

The three pairs with the highest absolute entries in the fitted
conditional concentration matrix are identical (AT\&T---Verizon,
Schlumberger---Baker Hughes  and Merrill Lynch---Morgan Stanley) for
the two methods and the 15 pairs with highest absolute entries in the
off-diagonal concentration matrix have an overlap of size 12. The
resulting graphs are slightly different although they share many
features.  Our graph has 136 edges, one more than that in the procedure
described in the paper, and   77 of the edges are shared.  Our graph
has more isolated vertices (15 vs. 9), slightly fewer cliques (62 vs.
81)  and the largest clique in our graph has six variables rather than
four. The graph is displayed to the left in Figure~\ref{fig:stable} and
features some clearly identified clusters of variables.

The selected graph is very unstable under bootstrap simulations. In the
spirit of \citet{meinshausen10stability}, we fit the graph on 2000
bootstrap samples. Only 28 edges are selected in more than half of
these samples. The resulting graph is shown in Figure~\ref{fig:stable}.
As many as 25 of these edges appear also as edges of the  estimator
proposed in (\ref{eq:orig}). It would have been interesting to be able
to compare with the same ``stability graph'' of the proposed procedure
but we suspect that they will match closely.

\section*{Latent directed structures}

In a sense the procedure described in this paper can be seen as a
modification of,\vadjust{\goodbreak} or an alternative to, factor analysis, in which
independent latent variables are sought to explain all the
correlations, corresponding to the graph for the observed variables
being completely empty.

Methods for identifying such models can, for example, be developed
using tetrad constraints [\citepp{Spirtes1993},
Drton, Sturmfels and Sullivant (\citeyear{drton2007algebraic})]. Another generalization of factor analysis
is to look for sparse \emph{directed} graphical models, which have now
been rather well establish\-ed through, for example, the FCI algorithm
[\citepp{Spirtes1993}, \citepp{Richardson2002}] with an algebraic
underpinning in \citeauthor{Sullivant2008} (\citeyear{Sullivant2008}).
Again this could be an alternative to the procedure described in this
interesting paper.

\section*{Summary}

We  effectively replaced the nuclear norm penalization of $L$ in the
paper by a fixed constraint on the rank. This might be easier to do
than choosing a reasonable value for the penalty on the trace of $L$.
Using this formulation, we could combine the EM algorithm with the
graphical lasso, enabling us to compute the solution with readily
available software. It would be interesting to see whether our
procedure can be shown to recover the correct sparsity structure under
similar assumptions to those in the paper. We want to congratulate the
authors again for a very interesting discussion paper.


\printaddresses

\end{document}